\documentclass{article}
\usepackage[english]{babel}
\usepackage[utf8]{inputenc}
\usepackage{geometry,amsmath,amssymb,enumerate,bbm,latexsym,theorem}
\usepackage{tikz}
\usepackage{hyperref}
\hypersetup
{
bookmarks=true,         
unicode=true,          
pdftoolbar=true,        
pdfmenubar=true,        
pdffitwindow=false,     
pdfstartview={FitH},    
pdftitle={},    
pdfauthor={Ali Abdallah},     
pdfsubject={},   
pdfcreator={Ali Abdallah},   
pdfproducer={Ali Abdallah}, 
pdfkeywords={} {} {} {}, 
pdfnewwindow=true,      
colorlinks=true,
linkcolor=black,
citecolor=red,
}

\newcommand{\nobracket}{}
\newcommand{\nocomma}{}
\newcommand{\tmSep}{; }
\newcommand{\tmaffiliation}[1]{\thanks{\textit{Affiliation:} #1}}
\newcommand{\tmop}[1]{\ensuremath{\operatorname{#1}}}
\newcommand{\tmstrong}[1]{\textbf{#1}}
\newcommand{\tmtextit}[1]{{\itshape{#1}}}
\newcommand{\tmverbatim}[1]{{\ttfamily{#1}}}
\newenvironment{enumeratenumeric}{\begin{enumerate}[1.] }{\end{enumerate}}
\newenvironment{enumerateroman}{\begin{enumerate}[i.] }{\end{enumerate}}
\newenvironment{itemizedot}{\begin{itemize} }{\end{itemize}}
\newenvironment{proof}{\noindent\textbf{Proof\ }}{\hspace*{\fill}$\Box$\medskip}
\newtheorem{corollary}{Corollary}
\newtheorem{definition}{Definition}
{\theorembodyfont{\rmfamily}\newtheorem{example}{Example}}
\newtheorem{proposition}{Proposition}
\newtheorem{theorem}{Theorem}

\begin{document}

\title{PUBLIC KEY CRYPTOGRAPHY BASED ON SOME EXTENSIONS OF GROUP}

\author{
  Ali Abdallah
  \tmaffiliation{Department of Mathematics{\tmSep}
  University of Rome Tor Vergata}
}

\date{March 14, 2016}

\maketitle

\begin{abstract}
  Bogopolski, Martino and Ventura in {\cite{bogopolski2010orbit}} introduced a
  general criteria to construct groups extensions with unsolvable conjugacy
  problem using short exact sequences. We prove that such extensions have
  always solvable word problem. This makes the proposed construction a
  systematic way to obtain finitely presented groups with solvable word
  problem and unsolvable conjugacy problem. It is believed that such groups
  are important in cryptography. For this, and as an example, we provide an
  explicit construction of an extension of Thompson group $F$ and we propose
  it as a base for a public key cryptography protocol.
\end{abstract}

\section{Introduction}

In 1997, Shor in his influential paper {\cite{shor1997polynomial}} proposed a
theoretical quantum algorithm for integer factorization into prime numbers
that runs in polynomial time. This would, in theory, compromise the current
most used public key crypto systems implementations (RSA, ECC, ...). Since
then, it is believed that group based cryptography might be a solution in
order to provide more secure cryptographic implementation
{\cite{myasnikov2008group}}. It is believed also that one solution could be to
find group with solvable word problem in linear time and with another very
hard decision problem.

The following group decision problems were first introduced by Max Dehn in
1911, within the context of closed 2-manifolds
\begin{itemizedot}
  \item The {\tmstrong{word problem}}: $\tmop{WP} (G) = \{x \in \Omega^{\ast}
  \mid x \underset{G}{=} 1\}$.
  
  \item The {\tmstrong{conjugacy problem}}: $\tmop{CP} (G) = \{x, y \in
  \Omega^{\ast} \mid y \underset{G}{=} g^{- 1} x g \nocomma, g \in G\} .$
\end{itemizedot}
Dehn showed that the word and conjugacy problems for the fundamental group of
a closed orientable surface of genus $g \geqslant 2$ is recursively solvable.
Furthermore, he defined the so called {\tmstrong{Dehn presentation}} for
groups for which he gives explicit algorithms for solving the word and
conjugacy problems. This is not the case for all finitely generated/presented
groups.

\begin{theorem}[Novikov {\cite{novikov1954algorithmic}}, Boone
{\cite{boone1954certain}}]
  There exists a finitely presented group whose word problem is recursively
  unsolvable.
\end{theorem}

A group with solvable conjugacy problem obviously has solvable word problem,
but the converse in not true in general. For example, Fidman
{\cite{fridman1969relation}} in 1960 showed that groups constructed earlier by
Novikov in {\cite{novikov1958unsolvability}} have unsolvable conjugacy and
solvable word problem. In this article we will be interested by such groups,
but using the technique developed by Bogopolski, Martino and Ventura
{\cite{bogopolski2010orbit}}, which permits to construct group that under some
conditions will have solvable/unsolvable conjugacy problem. In section 2, we
show that these extensions have always solvable word problem. In section 3, we
construct an explicit presentation of an extension of Thompson group $F$ with
solvable word problem and unsolvable conjugacy problem. Finally, in section 4
we illustrute the application of such extensions to cryptography.

\section{Group Extensions with Solvable Word Problem and Unsolvable Conjugacy
Problem}

Given a short exact sequence of groups
\begin{equation}
  1 \longrightarrow F \overset{\alpha}{\longrightarrow} G
  \overset{\beta}{\longrightarrow} H \longrightarrow 1 \label{ses} .
\end{equation}
We require that images and pre-images via the morphisms $\alpha$ and $\beta$
to be computable. The group of interest is the group $G$, extension of the
group $F$. \ Under suitable assumptions on the groups $F$ and $H$, \ the group
$G$ will have solvable or unsolvable conjugacy problem. In addition, we will
prove that in both cases the group $G$ will have solvable word problem.
Consider the following decidability problems:
\begin{itemizedot}
  \item The {\tmstrong{orbit decidability problem}}. Given a subgroup $A
  \leqslant \tmop{Aut} (F),$ we set
  \[ \tmop{OB} (F) = \{ (x, y) \in F \times F \mid \varphi (y) = g^{- 1} x g,
     \varphi \in A \tmop{and} g \in F\} . \]
  \item The {\tmstrong{$\varphi$-twisted conjugacy problem}}: For $\varphi \in
  \tmop{Aut} (F)$, we set
  \[ \tmop{TCP}_{\varphi} (F) = \{ (x, y) \in F \times F \mid y = g^{- 1} x
     \varphi (g) \} . \]
  \item {\tmstrong{The twisted conjugacy problem}} $\tmop{TCP} (F)$: $F$ has
  solvable twisted conjugacy problem (TCP) if $\tmop{TCP}_{\varphi}$ is
  solvable for every $\varphi \in \tmop{Aut} (F)$, and unsolvable TCP
  otherwise.
\end{itemizedot}
In the short exact sequence of (\ref{ses}), $\alpha (F)$ is a normal subgroup
of $G$. Thus we can identify $F$ with its image $\alpha (F)$ in $G$. In
addition, every inner automorphism of $G$ ($\varphi_g : G \longrightarrow G$,
that maps $x \in G$ to $g^{- 1} x g$), restricts to an automorphism of $F$. We
then define the {\tmstrong{action subgroup}} as follows:
\[ A_G = \{ \varphi_g \mid g \in G \} \leqslant \tmop{Aut} (F) . \]
\begin{theorem}
  {\cite{bogopolski2010orbit}} \label{ThExtensions}Suppose we are given a
  computable short exact sequence
  \[ 1 \longrightarrow F \overset{\alpha}{\longrightarrow} G
     \overset{\beta}{\longrightarrow} H \longrightarrow 1 \]
  with the following requirements:
  \begin{enumerateroman}
    \item for every $\varphi \in A_G$, $F$ has solvable $\varphi$-twisted
    conjugacy problem.,
    
    \item H has solvable conjugacy problem,
    
    \item for every $h \in H$ such that $h \neq 1$, the subgroup $\langle h
    \rangle$ has a finite index in its centralizer \ $C_H (h)$. And a
    corresponding coset representatives of $\langle h \rangle$ can be
    computed.
  \end{enumerateroman}
  Then:
  
  \begin{center}
    The conjugacy problem for $G$ is solvable if and only if the action
    subgroup $A_G $ is orbit decidable.
  \end{center}
\end{theorem}

\begin{theorem}
  \label{Hilte}Extensions obtained under the conditions of Theorem
  \ref{ThExtensions} have solvable word problem.
  
  \begin{proof}
    When the extension group $G$ has solvable conjugacy problem, the
    solvability of the word problem follows immediately.
    
    Now consider the case when the extension group $G$ has unsolvable
    conjugacy problem and let $y, y' \in G$, we map them into $H$, if they are
    not equal in $H$, then they cannot be equal in $G$. Otherwise $\beta (y
    y'^{- 1}) =_H 1$, and there is $f \in F$, such that $\alpha (f) = y y'^{-
    1}$. Since $\alpha$ is an injective map, it follows that $y =_G y'$ if and
    only if $f =_F 1$, which we can decide since the word problem of $F$ is
    solvable.
  \end{proof}
\end{theorem}

Theorem \ref{ThExtensions} gives a systematic effective method to construct
extensions of the group $F$ with unsolvable/solvable conjugacy problem for
which the word problem is solvable.

We need to introduce a new decision problem which we are going to link to the
orbit decidability problem. Let $B$ be a group and let $A$ be a subgroup of
$B$. The {\tmstrong{membership problem}} for $A$ in $B$ is defined as follows:
\[ \tmop{MP} (A, B) = \left\{ b \in B \mid a \in A \nocomma \nocomma, b
   \underset{A}{=} a \nocomma \right\} \]
That is, given an element $b \in B$ decide whether or not it belongs to $A$.

\begin{theorem}
  \label{MhTh}Let $G = \langle \nobracket x_1, \cdots, x_n \mid R_1, \cdots,
  R_m \rangle \nobracket$ a finitely presented group with unsolvable word
  problem. Let $A$ be the following subgroup of $F_n \times F_n$
  \[ A = \left\{ (x, y) \in F_n \times F_n \mid x \underset{G}{=} y \right\}
     \leqslant F_n \times F_n \]
  Then the membership problem of $A$ in $F_n \times F_n$ is unsolvable.
\end{theorem}

The construction presented in the above theorem is known as Mihailova's
construction and the group $A$ is called the Mihailova subgroup of $F_n \times
F_n$ associated with the group $G$ {\cite{mikhailova1966occurrence}}.

Let $F$ be a group, the {\tmstrong{stabilizer}} of a subgroup $A \leqslant F$
is:
\[ \tmop{Stab} (A) = \{ \varphi \in \tmop{Aut} (F) \mid \varphi (a) = a,
   \forall a \in A \} \leqslant \tmop{Aut} (F) . \]
We denote now by $\tmop{Stab}^{\ast} (A) = \tmop{Stab} (A) . \tmop{Inn} (F)
\leqslant \tmop{Aut} (F)$, the {\tmstrong{conjugacy stabilizer}} of $A$, where
$\tmop{Inn} (F)$ denotes the group of inner automorphism of $F$. The notation
"." denotes that elements of $\tmop{Stab}^{\ast} (A)$ are composition of an
element of $\tmop{Stab}^{} (A)$ and of an inner automorphism of $F$.

\begin{proposition}[{\cite{bogopolski2010orbit}}]
  Given a group $F$ and two subgroups $A \leqslant B \leqslant \tmop{Aut} (F)$
  and an element $v \in F$ such that $B \cap \tmop{Stab}^{\ast} (\langle v
  \rangle) = \{ \tmop{id} \}$. If $A \leqslant \tmop{Aut} (F) $ is orbit
  decidable then $\tmop{MP} (A, B)$ is solvable.
\end{proposition}

\begin{corollary}
  \label{FreeExtension}Suppose we are given a finitely presented group with
  solvable TCP, and such that $F_n \times F_n$ embeds in $\tmop{Aut} (F)$. If
  for $v \in F$, $\tmop{Stab}^{\ast} (< v >) \cap (F_n \times F_n) = \{
  \tmop{id} \}$, then it is possible to construct a finitely presented group
  with unsolvable CP, but solvable WP.
\end{corollary}

\section{Extension of Thompson Group $F$}

The structure of the automorphism group of $F$ was described by Brin in
{\cite{brin1996chameleon}}. As it is suggested by Brin, an easy way to
understand the automorphisms of $F$ is to look at $F$ as a subgroup of a
larger group. For this we introduce the group $\tmop{PL}_2 (\mathbbm{R})$ of
piece-wise-linear homeomorphisms of the real line, with dyadic breakpoints and
power of 2 slopes; allowing this time the set of breakpoints to be infinite,
but countable. In order to see $F$ as a subgroup of $\tmop{PL}_2
(\mathbbm{R})$, we conjugate elements of $F$ to the real line with a map
$\varphi : \mathbbm{R} \longrightarrow (0, 1)$ that is described in the
following figure:

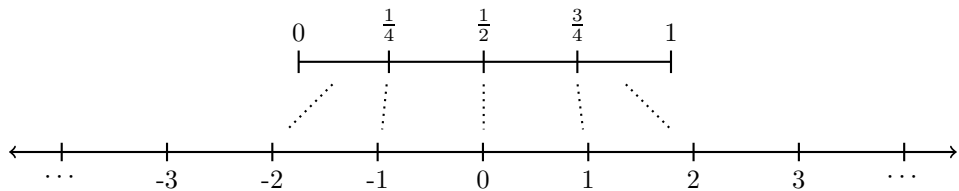
\begin{center}
\begin{tikzpicture}
\draw[-, thick,scale=3] (1.05,0.4) -- (2.7,0.4);
\draw[thick,scale=3] (1.45,0.35) -- (1.45,0.45) node[above] {$\frac{1}{4}$};
\draw[thick,scale=3] (1.87,0.35) -- (1.87,0.45) node[above] {$\frac{1}{2}$};
\draw[thick,scale=3] (1.05,0.35) -- (1.05,0.45) node[above] {0};
\draw[thick,scale=3] (2.285,0.35) -- (2.285,0.45) node[above] {$\frac{3}{4}$};
\draw[thick,scale=3] (2.7,0.35) -- (2.7,0.45) node[above] {1};

\draw[dotted, thick,scale=3] (1.87,0.30) -- (1.87,0.1);
\draw[dotted, thick,scale=3] (1.44,0.30) -- (1.42,0.1);
\draw[dotted, thick,scale=3] (2.286,0.30) -- (2.31,0.1);

\draw[dotted, thick,scale=3] (1.2,0.30) -- (1.0,0.1);
\draw[dotted, thick,scale=3] (2.5,0.30) -- (2.7,0.1);

\draw[<->, thick,scale=7] (-0.1,0) -- (1.7,0);
\foreach \x/\xtext in {0/$\cdots$,0.2/-3,0.4/-2,0.6/-1,0.8/0,1/1,1.2/2,1.4/3,1.6/$\cdots$}
    \draw[thick,scale=7] (\x,0.5pt) -- (\x,-0.5pt) node[below] {\xtext};
\end{tikzpicture}
\end{center}
 This defines an isomorphism from $F$ to a subgroup
of $\tmop{PL}_2 (\mathbbm{R})$ which is given by the conjugation map $f
\longmapsto \varphi^{- 1} f \varphi$. With this conjugation, it is easy to see
that elements of $F$ can be seen as elements of $\tmop{PL}_2 (\mathbbm{R})$,
with finitely many break points.

\begin{example}
  We illustrate in the following two diagrams the mapping of $A$ and $B$, the
  generators of Thompson group $F$. For instance the left diagram, which is is
  the mapping of $A$, can be read from top to bottom as $\varphi
  \longrightarrow A \longrightarrow \varphi^{- 1}$.
  
  
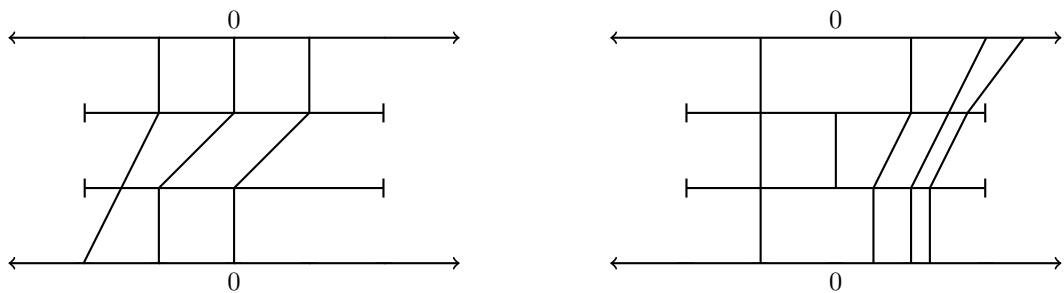
\begin{center}
\begin{tikzpicture}
\draw[<->, thick] (-3,0) -- (3,0);
\draw[thick] (0,0) -- (0,0) node {{\tiny |}}  node[above] {0};
\draw[thick] (-1,0) -- (-1,0) node {{\tiny |}}  node[above] {};
\draw[thick] (-2,0) -- (-2,0) node {{\tiny |}}  node[above] {};
\draw[thick] (1,0) -- (1,0) node {{\tiny |}}  node[above] {};
\draw[thick] (2,0) -- (2,0) node {{\tiny |}}  node[above] {};

\draw[|-|, thick] (-2,-1) -- (2,-1);
\draw[|-|, thick] (-2,-2) -- (2,-2);

\draw[<->, thick] (-3,-3) -- (3,-3);
\draw[thick] (0,-3) -- (0,-3) node {{\tiny |}}  node[below] {0};
\draw[thick] (-1,-3) -- (-1,-3) node {{\tiny |}}  node[below] {};
\draw[thick] (-2,-3) -- (-2,-3) node {{\tiny |}}  node[below] {};
\draw[thick] (1,-3) -- (1,-3) node {{\tiny |}}  node[below] {};
\draw[thick] (2,-3) -- (2,-3) node {{\tiny |}}  node[below] {};

\draw[thick] (-1,0) -- (-1,-1);
\draw[thick] (-1,-1) -- (-3/2,-2);
\draw[thick] (-3/2,-2) -- (-2,-3);

\draw[thick] (0,0) -- (0,-1);
\draw[thick] (0,-1) -- (-1,-2);

\draw[thick] (1,0) -- (1,-1);
\draw[thick] (1,-1) -- (0,-2);

\draw[thick] (-1,-2) -- (-1,-3);
\draw[thick] (0,-2) -- (0,-3);

\draw[<->, thick] (5,0) -- (11,0);
\draw (8,0) -- (8,0) node {{\tiny |}} node[above] {0};
\draw (7,0) -- (7,0) node {{\tiny |}} node[above] {};
\draw (6,0) -- (6,0) node {{\tiny |}} node[above] {};
\draw (9,0) -- (9,0) node {{\tiny |}} node[above] {};
\draw (10,0) -- (10,0) node {{\tiny |}} node[above] {};

\draw[|-|, thick] (6,-1) -- (10,-1);
\draw[|-|, thick] (6,-2) -- (10,-2);

\draw[<->, thick] (5,-3) -- (11,-3);
\draw (8,-3) -- (8,-3) node {{\tiny |}} node[below] {0};
\draw (7,-3) -- (7,-3) node {{\tiny |}} node[below] {};
\draw (6,-3) -- (6,-3) node {{\tiny |}} node[below] {};
\draw (9,-3) -- (9,-3) node {{\tiny |}} node[below] {};
\draw (10,-3) -- (10,-3) node {{\tiny |}} node[below] {};

\draw[thick] (7,0) -- (7,-1);
\draw[thick] (7,-1) -- (7,-2);
\draw[thick] (7,-2) -- (7,-3);

\draw[thick] (9,0) -- (9,-1);
\draw[thick] (8.5,-2) -- (8.5,-3);

\draw[thick] (8,-1) -- (8,-2);
\draw[thick] (9,-1) -- (8.5,-2);
\draw[thick] (9.5,-1) -- (9,-2);

\draw[thick] (10.5,0) -- (9.75,-1);
\draw[thick] (9.75,-1) -- (9.25,-2);
\draw[thick] (9.25,-2) -- (9.25,-3);

\draw[thick] (10,0) -- (9.5,-1);
\draw[thick] (9,-2) -- (9,-3);

\end{tikzpicture}
\end{center}

  \[ \begin{array}{lllllllllllll}
       \alpha (t) = t - 1 &  &  &  &  &  &  &  &  &  &  &  & \beta (t) =
       \left\{ \begin{array}{ll}
         t & t \leqslant 0\\
         t / 2 & 0 \leqslant t \leqslant 2\\
         t - 1 & t \geqslant 2
       \end{array} \right.
     \end{array}, \]
  and $\alpha, \beta$ generates a subgroup of $\tmop{PL}_2 (\mathbbm{R})$,
  which will be denoted by $\tmop{PL}_2 (I)$.
\end{example}

It follows that any element of the group $g = \varphi^{- 1} f \varphi \in
\tmop{PL}_2 (\mathbbm{R})$, must satisfy the following: $\exists M, N \in
\mathbbm{R}$ and $k, l \in \mathbbm{N}$ such that
\begin{itemize}
  \item For all $x > M$, $g (x) = x + k,$
  
  \item For all $x < N$, $g (x) = x + l.$
\end{itemize}
The following theorem is the key point to understand $\tmop{Aut} (F) \cong
\tmop{Aut} (\tmop{PL}_2 (I))$. A complete proof can be found in
{\cite{brin1996chameleon}}.

\begin{theorem}[Brin]
  Given $G \leqslant \tmop{PL}_2 (\mathbbm{R})$ we have:
  \[ \tmop{Aut} (G) \cong N (G), \]
  where $N (G)$ is the normalizer of $G$ in $\tmop{PL}_2 (\mathbbm{R}) .$
\end{theorem}

Viewing $F$ as subgroup of $\tmop{PL}_2 (\mathbbm{R})$, the automorphisms of
$F$ are elements of $\tmop{PL}_2 (\mathbbm{R})$ that conjugate $F$ to itself.
Let $\alpha \in \tmop{Aut} (G) \leqslant \tmop{PL}_2 (\mathbbm{R})$ be a
conjugator for $f_A$. There exists $g \in \tmop{PL}_2 (I)$, such that $f_A
\alpha = \alpha g$. Let $M$ be the bounded interval for $g$, such that $x > M
\nocomma$, $g (x) = x + l$ for some integer $l$. Take $x > M$, we have $f_A
(\alpha (x)) = \alpha (x) - 1 = \alpha g (x) = \alpha (x + l)$. By simply
writing down the general equation for $\alpha$, we can conclude that $l = -
1$. It follows that any $\alpha \in \tmop{Aut} (G)$, must satisfy $\alpha (x +
1) = \alpha (x) + 1$ outside some bounded interval.

In {\cite{burillo2013conjugacy}} the authors proved that the twisted conjugacy
problem is solvable for Thompson group $F$. They prove the existence of
extensions of Thompson group $F$ with unsolvable conjugacy problem. We do the
proof for the extensions construction parts in a slightly different way,
because it is helpful in constructing an explicit presentation of $F$.

\begin{theorem}[Burillo, Matucci, Ventura {\cite{burillo2013conjugacy}}]
  The twisted conjugacy problem is recursive for Thompson's group $F$.
\end{theorem}

The group $\tmop{PL}_2 (\mathbbm{R})$ contains an index two subgroup, namely
the subgroup of orientation preserving maps, usually denoted by $\tmop{PL}_2^+
(\mathbbm{R})$ and the subgroup of orientation reversing maps, usually denoted
by $\tmop{PL}_2^- (\mathbbm{R})$. In the same way, $\tmop{Aut}^+ (F)$ denotes
the subgroup of automorphisms of $F$ that preserve the orientation. This is
important as we are going to work only with orientation preserving
automorphism to obtain a short exact sequence between $F$, $\tmop{Aut}^+ (F)$
and $T \times T$. We recall that Thompson group $T$ is the group of piece-wise
linear homeomorphisms of \ $S^1 = [0, 1] / \{ 0 = 1 \}$ which are
differentiable with derivatives equal to powers of 2, except on a finite set
of dyadic rational numbers of the form $p / 2^q$.

Now by regarding $S^1$ as $\mathbbm{R}/\mathbbm{Z}$ we can view Thompson
group $T$ as $\tmop{PL}_2 (\mathbbm{R}/\mathbbm{Z})$. For an element $a \in
\tmop{Aut}^+ (F)$ we have $a (x + 1) = a (x) + 1$ outside of some bounded
interval $[M, N]$, we can map $a \mid_{(k - 1, k]}$ for some $k < M$ to the
quotient modulo $\mathbbm{Z}$ to obtain an element of $T$. In the same way, we
map $a \mid_{[l, l + 1)}$ to an element of $T$ for some $l > N$. We denote
this mapping by $\beta (a) = (a_-, a_+)$.

\begin{proposition}
  There is a short exact sequence:
  \[ 1 \longrightarrow F \overset{}{} \overset{i}{\longrightarrow}
     \tmop{Aut}^+ (F)  \overset{\beta}{\longrightarrow} T \times T
     \longrightarrow 1, \]
  where $i$ is the inclusion map.
  
  \begin{proof}
    Identifying $F$ with its image $i (F)$ in $\tmop{Aut}^+ (F)$, we can see
    that an element $f \in i (F)$ ($\exists M ; x > M$, $f (x) = x + k,$and
    $\exists N ; x < N$, $f (x) = x + l$) gets mapped to the identity by
    $\beta$, these elements are exactly $\ker (\beta)$.
    
    We shall see that $\beta$ is surjective. Fix $p < q \in \mathbbm{R}$, for
    $t = (t_-, t_+) \in T \times T$, and let $\widetilde{t_-}$ and
    $\widetilde{t_+}$ be a periodic lifting of $t_-$ and $t_+$ respectively
    such that $\widetilde{t_-} (p - 1) < p$ and $q < \widetilde{t_+} (q + 1)$.
    Next we compute $g_-, g, g_+ \in F$ such that $g_- (p - 1) =
    \widetilde{t_-} (p - 1)$, $g_- (p) = p, g (p) = p, g (q) = q$ and $g_+ (q)
    = q, g_+ (q + 1) = \widetilde{t_+} (x)$.
    \[ a (x) = \left\{ \begin{array}{ll}
         \widetilde{t_-} (x) & x \leqslant p - 1\\
         g_- (x) & p - 1 \leqslant x \leqslant p\\
         g (x) & p \leqslant x \leqslant q\\
         g_+ (x) & q \leqslant x \leqslant q + 1\\
         \widetilde{t_+} (x) & q + 1 \leqslant x,
       \end{array} \right. \]
    It is clear that $a \in \tmop{Aut}^+ (F)$ and $\beta (a) = t = (t_-, t_+)
    .$
  \end{proof}
\end{proposition}

\begin{theorem}
  \label{ThExt}For every free subgroup $F_2 \leqslant T$ of rank $2$, there is
  an extension of Thompson group $F$ with solvable word problem and unsolvable
  conjugacy problem.
  
  \begin{proof}
    Let $F_2 \cong \langle u, v \rangle \leqslant T$ a free group of rank 2 in
    $T$. We construct the following free product $F_2 \times F_2 \cong \langle
    a, b \rangle \times \langle c, d \rangle \leqslant T \times T$, where $a =
    u^2, b = v^2, c = u v u^{- 1}$ and $d = v u v^{- 1}$.
    
    We construct lifts $\hat{a}, \hat{b}, \hat{c}, \hat{d}$ of $a, b, c$ and
    $d$ respectively such that $\beta (\hat{a}) = (a, 1), \beta (\hat{b}) =
    (b, 1), \beta (\hat{c}) = (1, c)$ and $\beta (\hat{d}) = (1, d)$. By
    construction this gives a copy $B : = F_2 \times F_2 \cong \langle
    \hat{a}, \hat{b}, \hat{c}, \hat{d} \rangle \leqslant \tmop{Aut}^+ (F)$.
    
    Now in order to conclude the proof, we need to show that for some $v \in
    F$ one has $B \cap \tmop{Stab}^{\ast} (v) = \{ \tmop{id} \}$. For that let
    $v \in F, v (x) = x + 1, \forall x \in \mathbbm{R}$ and let $s \in B \cap
    \tmop{Stab}^{\ast} (v)^{}$. Since $s$ acts on $F$ by conjugation, we have
    $s^{- 1} v s = g^{- 1} v g$ for some $g \in F$, thus $s g^{- 1}$ and $v$
    commutes. Since $v$ is periodic of period $1$ over the entire real line,
    it follows that $\beta (s g^{- 1}) = (t, t) = \beta (s)$, for some $t \in
    T$. By writing down $s$ in terms of reduced words over $s g^{- 1} = w_1
    (\hat{a}, \hat{b}) w_2 (\hat{c}, \hat{d}) \in \langle \hat{a}, \hat{b},
    \hat{c}, \hat{d} \rangle$, we can conclude that $\beta (s) = (t, t)$ can
    only happen if and only if $t = 1_T$. And so $B \cap \tmop{Stab}^{\ast}
    (v) = \{ \tmop{id} \}$. Therefore we can apply Theorem \ref{FreeExtension}
    to obtain an extension with solvable word problem and unsolvable conjugacy
    problem.
  \end{proof}
\end{theorem}

It is known that Thompson group $F$ does not contain a subgroup which is
isomorphic to the free group of rank 2. Unlike Thompson group $F$, the group
$T$ contains free subgroups of rank 2. This implies the existence of
extensions of $F$ with unsolvable conjugacy problem.

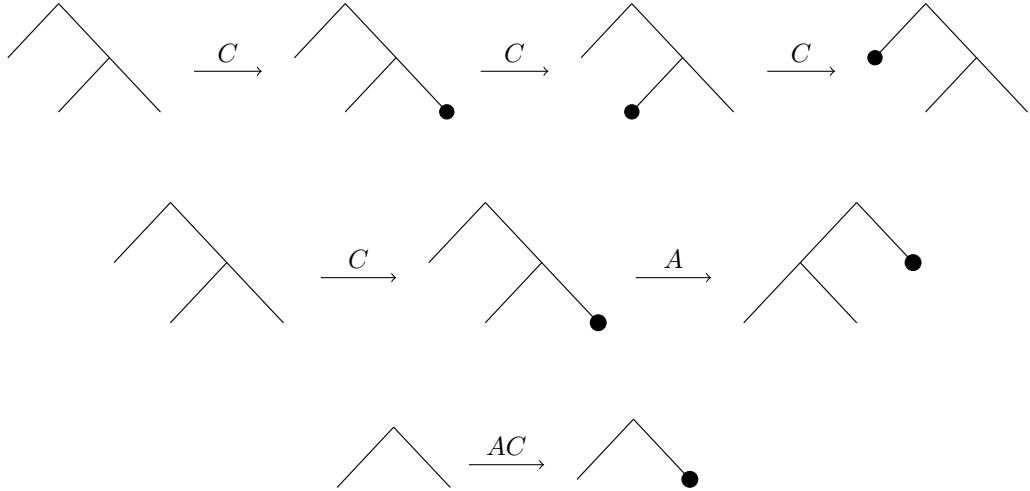
\begin{center}
\begin{tabular}{l c c r}
\begin{tikzpicture}[level distance=8mm,scale=0.9]
\coordinate
child {[fill] circle (0pt)}
child {[fill] circle (0pt) child {[fill] circle (0pt)} child {[fill] circle (0pt) child[fill=none] {edge from parent[draw=none]}}};
\draw[->] (2cm,-1cm)  --  (3cm,-1cm) node[midway,sloped,above] {$C$} ;
\end{tikzpicture} &

\begin{tikzpicture}[level distance=8mm,scale=0.9]
\coordinate
child {[fill] circle (0pt)}
child {[fill] circle (0pt) child {[fill] circle (0pt)} child {[fill] circle (3pt) child[fill=none] {edge from parent[draw=none]}}};
\draw[->] (2cm,-1cm)  --  (3cm,-1cm) node[midway,sloped,above] {$C$} ;
\end{tikzpicture} &

\begin{tikzpicture}[level distance=8mm,scale=0.9]
\coordinate
child {[fill] circle (0pt)}
child {[fill] circle (0pt) child {[fill] circle (3pt)} child {[fill] circle (0pt) child[fill=none] {edge from parent[draw=none]}}};
\draw[->] (2cm,-1cm)  --  (3cm,-1cm) node[midway,sloped,above] {$C$} ;
\end{tikzpicture} &

\begin{tikzpicture}[level distance=8mm,scale=0.9]
\coordinate
child {[fill] circle (3pt)}
child {[fill] circle (0pt) child {[fill] circle (0pt)} child {[fill] circle (0pt) child[fill=none] {edge from parent[draw=none]}}};
\end{tikzpicture} 
\end{tabular}
\label{fig:t-c3-gen}
\end{center}

\begin{center}
\begin{tabular}{l c r }
\begin{tikzpicture}[level distance=8mm]
\coordinate
child {[fill] circle (0pt)}
child {[fill] circle (0pt) child {[fill] circle (0pt)} child {[fill] circle (0pt) child[fill=none] {edge from parent[draw=none]}}};
\draw[->] (2cm,-1cm)  --  (3cm,-1cm) node[midway,sloped,above] {$C$} ;
\end{tikzpicture} &

\begin{tikzpicture}[level distance=8mm]
\coordinate
child {[fill] circle (0pt)}
child {[fill] circle (0pt) child {[fill] circle (0pt)} child {[fill] circle (3pt) child[fill=none] {edge from parent[draw=none]}}};
\draw[->] (2cm,-1cm)  --  (3cm,-1cm) node[midway,sloped,above] {$A$} ;
\end{tikzpicture} &

\begin{tikzpicture}[level distance=8mm]
\coordinate
child {[fill] circle (0pt) child {[fill] circle (0pt)} child {[fill] circle (0pt) child[fill=none] {edge from parent[draw=none]}}}
child {[fill] circle (3pt)};
\end{tikzpicture} 
\label{fig:t-ac-gen}
\end{tabular}
\end{center}

\begin{center}
\begin{tabular}{l r }
\begin{tikzpicture}[level distance=8mm]
\coordinate
child {[fill] circle (0pt)}
child {[fill] circle (0pt)};
\draw[->] (1cm,-0.5cm)  --  (2cm,-0.5cm) node[midway,sloped,above] {$AC$} ;
\end{tikzpicture} &
\begin{tikzpicture}[level distance=8mm]
\coordinate
child {[fill] circle (0pt)}
child {[fill] circle (3pt)};
\end{tikzpicture}
\label{fig:t-ac2-gen}
\end{tabular}
\end{center}

From the computation presented in the above diagrams, the following equations
hold in $T$:
\begin{enumerateroman}
  \item $C^3 = 1,$
  
  \item $(A C)^2 = 1.$
\end{enumerateroman}
In order to find a copy of $F_2$ in $T$, we view Thompson group $T$ as the
group of orientation preserving homeomorphisms of the real projective line
$\mathbbm{R}P^1$, which are piecewise $\tmop{PSL}_2 (\mathbbm{Z})$ and
differentiable except on a finite set of rational numbers. In this way, the
group $\tmop{PSL}_2 (\mathbbm{Z})$ can be seen as a subgroup of $T$. The
standard presentation of $\tmop{PSL}_2 (\mathbbm{Z})$ is:

\begin{center}
  $\tmop{PSL}_2 (\mathbbm{Z}) = \langle a, b \mid a^2 = b^3 = 1 \rangle$
\end{center}

\begin{proposition}
  {\cite{fossas2014thompson}} The subgroup $H = \langle a b a b^{- 1}, a b^{-
  1} a b \rangle$ is a free non abelian group of rank 2 in $T$.
\end{proposition}

We can explicitly obtain an isomorphic copy of $H$ in $T$ using what we have
computed before ($C^3 = (A C)^2 = 1$). Thus the following copy of $H$ in $T$
written as follows
\[ H = \langle u = A C^2 A, v = A^2 C^2  \rangle \]
is a free non abelian group of rank 2 in $T$. We take $\langle a, b \rangle
\times \langle c, d \rangle \leqslant T \times T$, where $a = u^2, b = v^2, c
= u v u^{- 1}$ and $d = v u v^{- 1}$. And so we have a copy of $F_2 \times F_2
\cong \langle \hat{a}, \hat{b} \rangle \times \langle \hat{c}, \hat{d} \rangle
\leqslant \tmop{Aut}^+ (F)$. For a product of the free group of rank 2 on the
same generators $F_2 = \langle \nobracket x, y \rangle \nobracket$, we have an
embedding of $F_2 \times F_2$ into $\tmop{Aut}^+ (F)$.
\[ \begin{array}{llll}
     & F_2 \times F_2 : & \longrightarrow & \tmop{Aut}^+ (F) \leqslant
     \tmop{Aut} (F)\\
     & (w_1, w_2) & \longmapsto & \widehat{w_1} \widehat{w_2}
   \end{array}, \]
Let $G = \langle x, y \mid R_1, \cdots, R_m \rangle$ be a finitely presented
group on two generators with unsolvable word problem (for example the group
presented in {\cite{wanggroups}}). We then construct a group
\[ A = \left\{ (\phi, \psi) \in F_2 \times F_2 \leqslant \tmop{Aut}^+ (F)
   \mid \phi \underset{G}{=} \psi \right\} \leqslant F_2 \times F_2 \]
As in Theorem \ref{MhTh}, the group $A$ is our action sub-group of $\tmop{Aut}
(F)$ with unsolvable orbit decidability problem (since MP$(A, F_2 \times F_2)$
is unsolvable) and it is finitely generated:
\[ A \cong \langle \phi_1 = (1, R_1), \cdots, \phi_m = (1, R_m), \phi_{m + 1}
   = (x, x), \phi_{m + 2} = (y, y) \rangle \]
These generators get mapped into the copy of $F_2 \times F_2$ in $\tmop{Aut}^+
(F)$:
\[ A \cong \langle \widehat{\phi_1} = \widehat{R_1}, \cdots, \widehat{\phi_m}
   = \widehat{R_m}, \widehat{\phi_{m + 1}} = \hat{x} \hat{x}, \widehat{\phi_{m
   + 2}} = \hat{y} \hat{y} \rangle \]
Let $F_n = \langle t_1, \cdots, t_n \rangle$ be the free group of rank $n$
(free groups have solvable conjugacy problem and $[C_{F_n} : a] = 1$, for $a
\in F_n$). The group $G$ given by the following presentation:
\begin{equation}
  G = \left\langle \alpha, \beta, t_1, \cdots, t_n  \middle|  \begin{array}{l}
    {}[\alpha \beta^{- 1}, \alpha^{- 1} \beta \alpha], [\alpha \beta^{- 1},
    \alpha^{- 2} \beta \alpha^2],\\
    t_j^{- 1} \alpha t_j = \widehat{\phi_j}^{- 1} \alpha \widehat{\phi_j},\\
    t_j^{- 1} \beta t_j = \widehat{\phi_j}^{- 1} \beta \widehat{\phi_j}
  \end{array} \right\rangle, \label{ExtT}
\end{equation}
$j = 1, \cdots, n$, is an extension of Thompson group $F$ with solvable word
problem and unsolvable conjugacy problem.

\section{Application to Cryptography}

\begin{definition}
  Let $G$ be a group with a finite generating set $S$. The {\tmstrong{growth
  function}} $\gamma (n)$ is defined for every $n \in \mathbbm{N}$ as the
  number of elements of $G$ which are product of at most $n$ elements of $S$.
\end{definition}

The $\lim_{n \rightarrow \infty} \sqrt[n]{\gamma (n)}$ exists always
{\cite{de2000topics}}. A finitely generated group is said to have
{\tmstrong{exponential growth}} if the limit is positive, and
{\tmstrong{subexponential growth}} if the limit is $0$. Exponential growth
property is very important when using groups as base for cryptographic
algorithm implementation, because groups with exponential growth provide
larger key space. It is also interesting to note that having exponential or
sub-exponential growth is an intrinsic property of the group, that is, it does
not depend upon the finite generating set {\cite{de2000topics}}.

\begin{proposition}
  Thompson group $F$ has exponential growth.
  
  \begin{proof}
    Consider words of $F$ in the following form:
    \[ X_1^{\epsilon_1} X_0^{- 1} X_1^{\epsilon_2} X_0^{- 1} \cdots X_0^{- 1}
       X_1^{\epsilon_n}, \]
    where $\epsilon_1, \cdots, \epsilon_n \in \mathbbm{N}$. By using the
    infinite presentation of $F$ we can move the $X_0^{- 1}$ to the right to
    obtain:
    \begin{equation}
      X_1^{\epsilon_1} X_2^{\epsilon_2} \cdots X_n^{\epsilon_n} X_0^{- (n -
      1)}
    \end{equation}
    Words of $(5.1)$ are in their normal forms, therefore we get a different
    element of $F$ for every different values of $\epsilon_1, \cdots,
    \epsilon_n$.
  \end{proof}
\end{proposition}

Since Thompson group $F$ is contained in its extensions, we have:

\begin{corollary}
  Any extension of Thompson group $F$ has exponential growth.
\end{corollary}

We have seen that the extensions of Thompson group $F$ developed in the
previous chapter have all the interesting properties to be a base of a crypto
system (exponential growth, solvable word problem, unsolvable conjugacy
problem). In this section we discuss a generic application of such extensions
for public-key cryptography. The protocol that we are going to describe is due
to Anshel, Anshel and Goldfeld {\cite{anshel1999algebraic}}. The only
necessary requirement for the algorithm is the solvability of the word
problem. The security of the protocol relies on the so called
{\tmstrong{simultaneous conjugacy problem}}, which is harder than the
conjugacy problem. It can be stated as follows:
\[ \tmop{SCP} (G) = \{ (x_1, \cdots, x_n, y_1, \cdots, y_n) \in G^{2 n} \mid
   \exists g \in G ; y_i = g^{- 1} x_i g \} . \]
The following proposition is immediate:

\begin{proposition}
  For a group $G$ we have the following reduction:
  \[ \tmop{CP} (G) \leqslant_m \tmop{SCP} (G) \]
\end{proposition}

and so as a corollary:

\begin{corollary}
  The simultaneous conjugacy problem is unsolvable for extensions of Thompson
  group $F$ with unsolvable conjugacy problem.
\end{corollary}

Given $s_1, \cdots, s_m$ and $t_1, \cdots, t_n$ are elements of a group $G$.
The protocol can be described as follows between two users \tmverbatim{A} and
\tmverbatim{B}.
\begin{enumeratenumeric}
  \item User \tmverbatim{A} computes a random element $a \in G$ as a word in
  terms of $s_1, \cdots, s_m$ denoted by $a (s_1, \cdots, s_m)$ and sends
  $a^{- 1} t_1 a, \cdots, a^{- 1} t_n a$ to \tmverbatim{B}.
  
  \item User \tmverbatim{B} computes a random element $b \in G$ as a word in
  terms of $t_1, \cdots, t_n$ denoted by $b (t_1, \cdots, t_n)$ and sends
  $b^{- 1} s_1 b, \cdots, b^{- 1} s_m b$ to \tmverbatim{A}.
  
  \item \tmverbatim{A} computes $a (b^{- 1} s_1 b, \cdots, b^{- 1} s_m b) =
  b^{- 1} a b$.
  
  \item \tmverbatim{B} computes $b (a^{- 1} t_1 a, \cdots, a^{- 1} t_n a) =
  a^{- 1} b a$.
\end{enumeratenumeric}
At this point, user \tmverbatim{A} disposes of the private key $a^{- 1} b^{-
1} a b = [a, b]$ and user \tmverbatim{B} the key $b^{- 1} a^{- 1} b a = [b, a]
= [a, b]^{- 1}$ and user B can easily compute $[a, b]$ sharing this way the
same private secret with \tmverbatim{A}. An adversary \tmverbatim{C} observing
the transmissions $(1)$ and $(2)$ is unable to figure out $a$ and $b$ unless
he/she can solve the set of simultaneous conjugacy over $G$.

One may think of applying the same procedure to create extensions of $B_n$
with solvable word problem and unsolvable conjugacy problem which makes a
cryptographic protocol based on extensions of $B_n$ more difficult to break.
Recently Meneses and Ventura in {\cite{gonzalez2014twisted}} proved that:

\begin{theorem}
  The twisted conjugacy problem is solvable for $B_n$.
\end{theorem}

Dyer and Grossman extensively studied the automorphisms of $B_n$:

\begin{theorem}[Dyer, Grossman {\cite{dyer1981automorphism}}]
  $\tmop{Aut} (B_n) = \tmop{Inn} (B_n) \sqcup \tmop{Inn} (B_n) . \epsilon$,
  where $\epsilon : B_n \longrightarrow B_n$ is the automorphism which inverts
  the generators of $B_n$ ($\sigma_i \longrightarrow \sigma_i^{- 1}$).
\end{theorem}

This means that given $\varphi \in \tmop{Aut} (B_n)$ and by the previous
theorem, either $\varphi = \gamma_g, g \in B_n$ or $\varphi = \gamma_g .
\epsilon, g \in B_n \nocomma$, where $\gamma_g$ is the conjugation map. Given
a finitely generated subgroup $A = \langle \varphi_1, \cdots, \varphi_m
\rangle \leqslant \tmop{Aut} (B_n)$, for every $i = 1, \cdots, m$ we have
$\varphi_i = \gamma_{g_i} \epsilon^{\alpha}$ ($\alpha = 0$ or $1$). With this,
it is easy to see that for $u, v \in B_n$, deciding whether or not there
$\varphi \in A$ and $g \in B_n$ such that $\varphi (u) = g^{- 1} v g$ reduces
to solving conjugacy problem in $B_n$, which is solvable:

\begin{corollary}
  Every finitely generated subgroup $A \leqslant \tmop{Aut} (B_n)$ is orbit
  decidable.
\end{corollary}

Thus what we have applied on Thompson group $F$ to obtain extensions with
unsolvable conjugacy problem does not apply on the braid group $B_n$ and we
have:

\begin{corollary}
  All extensions of $B_n$ that can be constructed under the conditions of
  Theorem \ref{ThExtensions} have solvable conjugacy problem.
\end{corollary}

\


\begin{thebibliography}{WXLL14}
  \bibitem[AAG99]{anshel1999algebraic}Iris Anshel, Michael Anshel, and Dorian
  Goldfeld. {\newblock}An algebraic method for public-key cryptography.
  {\newblock}\tmtextit{Mathematical Research Letters}, 6:287--292, 1999.
  
  \bibitem[BMV10]{bogopolski2010orbit}Oleg Bogopolski, Armando Martino, and
  Enric Ventura. {\newblock}Orbit decidability and the conjugacy problem for
  some extensions of groups. {\newblock}\tmtextit{Transactions of the American
  Mathematical Society}, 362(4):2003--2036, 2010.
  
  \bibitem[BMV13]{burillo2013conjugacy}José Burillo, Francesco Matucci, and
  Enric Ventura. {\newblock}The conjugacy problem in extensions of Thompson's
  group F. {\newblock}\tmtextit{arXiv preprint arXiv:1307.6750}, 2013.
  
  \bibitem[Boo54]{boone1954certain}William~W Boone. {\newblock}Certain simple,
  unsolvable problems of group theory. i. {\newblock}In \tmtextit{Indagationes
  Mathematicae (Proceedings)}, volume~57, pages 231--237. Elsevier, 1954.
  
  \bibitem[Bri96]{brin1996chameleon}Matthew~G Brin. {\newblock}The chameleon
  groups of Richards J. Thompson: automorphisms and dynamics.
  {\newblock}\tmtextit{Publications Mathématiques de l'IHÉS}, 84:5--33,
  1996.
  
  \bibitem[DG81]{dyer1981automorphism}Joan~L Dyer and Edna~K Grossman.
  {\newblock}The automorphism groups of the braid groups.
  {\newblock}\tmtextit{American Journal of Mathematics}, pages 1151--1169,
  1981.
  
  \bibitem[dLH00]{de2000topics}Pierre de~La~Harpe. {\newblock}\tmtextit{Topics
  in geometric group theory}. {\newblock}University of Chicago Press, 2000.
  
  \bibitem[Fos14]{fossas2014thompson}Ariadna Fossas. {\newblock}Thompson's
  group t, undistorted free groups and automorphisms of the flip graph.
  {\newblock}In \tmtextit{Extended Abstracts Fall 2012}, pages 45--49.
  Springer, 2014.
  
  \bibitem[Fri69]{fridman1969relation}AA~Fridman. {\newblock}On the relation
  between the word problem and the conjugacy problem in finitely defined
  groups. {\newblock}1969.
  
  \bibitem[GMV14]{gonzalez2014twisted}Juan González-Meneses and Enric
  Ventura. {\newblock}Twisted conjugacy in braid groups.
  {\newblock}\tmtextit{Israel Journal of Mathematics}, 201(1):455--476, 2014.
  
  \bibitem[Mik66]{mikhailova1966occurrence}KA~Mikhailova. {\newblock}The
  occurrence problem for direct products of groups.
  {\newblock}\tmtextit{Matematicheskii Sbornik}, 112(2):241--251, 1966.
  
  \bibitem[MSU08]{myasnikov2008group}Alexei Myasnikov, Vladimir Shpilrain, and
  Alexander Ushakov. {\newblock}\tmtextit{Group-based cryptography}.
  {\newblock}Springer Science \& Business Media, 2008.
  
  \bibitem[Nov54]{novikov1954algorithmic}PS~Novikov. {\newblock}On algorithmic
  unsolvability of the word problem. {\newblock}1954.
  
  \bibitem[Nov58]{novikov1958unsolvability}Petr~Sergeevich Novikov.
  {\newblock}Unsolvability of the conjugacy problem in group theory.
  {\newblock}1958.
  
  \bibitem[Sho97]{shor1997polynomial}Peter~W Shor. {\newblock}Polynomial-time
  algorithms for prime factorization and discrete logarithms on a quantum
  computer. {\newblock}\tmtextit{SIAM journal on computing}, 26(5):1484--1509,
  1997.
  
  \bibitem[WXLL14]{wanggroups}Xiaofeng Wang, Chen Xu, Guo Li, and Hanling Lin.
  {\newblock}Groups with two generators having unsolvable word problem and
  presentations of Mihailova subgroups. {\newblock}Technical report,
  Cryptology ePrint Archive, Report 2014.528, http://eprint. iacr. org, 2014.
\end{thebibliography}
\end{document}